\documentclass[12pt]{amsart}
\usepackage{amsmath,amssymb,amsthm}

\newcommand{\ds}{\displaystyle}
\newcommand{\F}{\mathcal F}
\newcommand{\N}{\mathcal N}
\newcommand{\W}{\mathcal W}

\theoremstyle{definition}

\theoremstyle{remark}
\begin{document}

\title[ON THE CLASSIFICATION OF THE ALMOST CONTACT METRIC MANIFOLDS]
{ON THE CLASSIFICATION OF THE ALMOST CONTACT METRIC MANIFOLDS}%

\author{Valentin \, A. Alexiev and Georgi \, T. Ganchev}

\address{Bulgarian Academy of Sciences, Institute of Mathematics and Informatics,
Acad. G. Bonchev Str. bl. 8, 1113 Sofia, Bulgaria}%
\email{ganchev@math.bas.bg}%


\begin{abstract}

The vector space of the tensors $\mathcal F$ of type (0,3) having the same
symmetries as the covariant derivative of the fundamental form of an almost
contact metric manifold is considered. A scheme of decomposition of $\F$
into orthogonal components which are invariant under the action of $U(n)\times 1$
is given. Using this decomposition there are found 12 natural basic classes of
almost contact metric manifolds. The classes of cosymplectic, $\alpha$-Sasakian,
$\alpha$-Kenmotsu, etc. manifolds fit nicely to these considerations. On the
other hand, many new interesting classes of almost contact metric manifolds arise.
\end{abstract}

\maketitle

\thispagestyle{empty}

\section{Preliminaries}

Let $V$ be a $(2n+1)$-dimensional real vector space with almost contact metric
structure $(\varphi, \xi, \eta, g)$, where $\varphi$ is a tensor of type (1,1),
$\xi$ is a vector, $\eta$ is a covector and $g$ is a definite metric so that
$$\begin{array}{l}
\varphi^2x=-x+\eta(x)\,\xi, \quad \varphi(\xi)=0, \quad \eta \circ \varphi =0,\\
[2mm]
g(\xi, \xi)=1, \quad g(\varphi x, \varphi y)=g(x,y) - \eta(x) \eta(y)
\end{array}$$
for arbitrary vectors  $x$, $y$ in $V$. For arbitrary $x \in V$, we denote
$hx=-\varphi^2x$.

We consider the subspace $\F$ of $V^*\otimes V^*\otimes V^*$ defined
by the conditions:
$$(1) \; \F=\{F\in \F\,|\,F(x,y,z)=-F(x,z,y)
=-F(x,\varphi y, \varphi z)+\eta(y)F(x, \xi, z) +\eta(z)F(x, y, \xi)\}$$
for $x, y, z$ in $V$.

Let $\{ e_i\}, i=1,...,2n+1$ be an orthonormal basis of $V$. The metric $g$
induces an inner product in the vector space $\mathcal F$:
$$<F',F''>=\sum_{i,j,k=1}^{2n+1}F'(e_i,e_j,e_k)F''(e_i,e_j,e_k); \quad F',F''\in \F.$$

We associate with every $F \in \mathcal F$ the following covectors:
$$(2) \; f(F)(z)=\sum_{i}F(e_i,e_i,z), \quad f^*(F)(z)=\sum_{i}F(e_i,\varphi \,e_i,z),
\quad \omega(F)(z)=F(\xi,\xi,z).$$

The standard representation of $U(n)\times 1$ in $V$ induces an associated representation of
$U(n)\times 1$ in $\mathcal F$.

It is well known the following
\vskip 2mm
{\bf Lemma 1.}
Let $L$ be an involutive isometry of $\F$, which commutes with the action of
$U(n)\times 1$ in $\F$. Then
$$\F=L^+(\F)\oplus L^-(\F),$$
where $L^+(\F)$ and $L^-(\F)$ are the eigenspaces of $L$ corresponding to the eigenvalues
$+1$ and $-1$ of $L$. The decomposition is orthogonal and invariant under the action of
$U(n)\times 1$. The components of an element $F\in \F$ in $L^+(\F)$ and $L^-(\F)$ are
$$F^+=\frac{1}{2}(F+LF), \quad F^-=\frac{1}{2}(L-LF).$$

\section{Associated forms with an element of $\F$}

With every $F$ in $\F$ we associate the following basic forms:
$$\begin{array}{rl}
F_1(F)(x,y,z)=&\eta(x)F(\xi,y,z),\\
[2mm]
F_2(F)(x,y,z)=&\eta(y)F(x,\xi,z)-\eta(z)F(x,\xi,y),\\
[2mm]
F_3(F)(x,y,z)=&\eta(x)\eta(y)F(\xi,\xi,z)-\eta(x)\eta(z)F(\xi,\xi,y)\\
[2mm]
 =&\eta(x)\eta(y)\omega(F)(z)-\eta(x)\eta(z)\omega (F)(y),\\
[2mm]
hF(x,y,z)=&F(hx,hy,hz).
\end{array}$$
\vskip 2mm
{\bf Lemma 2.} Let $F \in \F$. Then $F_i(F) (i=1,2,3)$ and $hF$ are elements of $\F$ and
$$F=hF+F_1(F)+F_2(F)-F_3(F).$$

Further we consider the forms

$$\begin{array}{l}
F_4(F)(x,y,z)=\eta(y)F(\varphi x,\xi,\varphi z)-\eta(z)F(\varphi x, \xi,\varphi y),\\
[2mm]
F_5(F)(x,y,z)=\eta(y)F(z,\xi,x)-\eta(z)F(y,\xi,x),\\
[2mm]
F_6(F)(x,y,z)=\eta(y)F(\varphi z,\xi,\varphi x)-\eta(z)F(\varphi y,\xi,\varphi x),\\
[2mm]
F_7(F)(x,y,z)=\ds{\frac{1}{2n}f(F)(\xi)\{\eta(z)g(x,y)-\eta(y)g(x,z)\}},\\
[2mm]
F_8(F)(x,y,z)=-\ds{\frac{1}{2n}f^*(F)(\xi)\{\eta(z)g(x, \varphi y)-\eta(y)g(x,\varphi z)\}}
\end{array}$$
associated with an arbitrary $F\in \F$.
\vskip 2mm
{\bf Lemma 3.} Let $F\in \F$. Then $F_i(F)\,(i=4,5,6,7,8)$ are elements of $\F$.
\vskip 2mm
{\bf Lemma 4.} Let $F\in \F$. The following relations are valid
$$\begin{array}{l}
F_{11}(F)=F_1(F), F_{12}(F)=F_3(F), F_{13}(F)=F_3(F), F_{14}(F)=0, F_{15}(F)=0,\\
[2mm]
F_{21}(F)=F_3(F), F_{22}(F)=F_2(F), F_{23}(F)=F_3(F), F_{24}(F)=F_4(F), F_{25}(F)=F_5(F),\\
[2mm]
F_{31}(F)=F_3(F), F_{32}(F)=F_3(F), F_{33}(F)=F_3(F), F_{34}(F)=0, F_{35}(F)=0,\\
[2mm]
F_{41}(F)=0, F_{42}(F)=F_4(F), F_{43}(F)=0, F_{44}(F)=F_2(F)-F_3(F), F_{45}(F)=F_6(F),\\
[2mm]
F_{51}(F)=0, F_{52}(F)=F_5(F), F_{53}(F)=0, F_{54}(F)=F_6(F), F_{55}(F)=F_2(F)-F_3(F),\\
[2mm]
F_{71}(F)=F_{17}(F)=F_{73}(F)=F_{37}(F)=0, F_{7i}(F)=F_{i7}(F)=F_7(F), \, i=2,4,5,7,\\
[2mm]
F_{81}(F)=F_{18}(F)=F_{83}(F)=0, F_{58}(F)=F_{85}(F)=-F_{48}(F)=-F_{84}(F)=F_{88}(F)=F_8(F),\\
[2mm]
h(F_i(F)=F_i(hF)=0,\,i=1,...,8,
\end{array}$$
where $F_{ij}(F)=F_i(F_j(F))$.
\vskip 2mm
{\bf Lemma 5.} Let $F\in \F$. Then we have
$$\begin{array}{l}
f(F_1(F))=\omega(F), \; f^*(F_1(F))=0, \; \omega(F_1(F))=\omega(F),\\
[2mm]
f(F_2(F))=\omega(F)+f(F)(\xi)\eta, \; f^*(F_2(F))=f^*(F)(\xi)\eta, \; \omega(F_2(F))=\omega(F),\\
[2mm]
f(F_3(F))=\omega(F), \; f^*(F_3(F))=0, \; \omega(F_3(F))=\omega(F),\\
[2mm]
f(F_4(F))=f(F)(\xi)\eta, \; f^*(F_4(F))=f^*(F)(\xi)\eta, \; \omega(F_4(F))=0,\\
[2mm]
f(F_5(F))=f(F)(\xi)\eta, \; f^*(F_5(F))=-f^*(F)(\xi)\eta, \; \omega(F_5)(F)=0,\\
[2mm]
f(F_6(F))=f(F)(\xi)\eta, \; f^*(F_6(F))=-f^*(F)(\xi)\eta, \; \omega(F_6(F))=0,\\
[2mm]
f(F_7(F))=f(F)(\xi)\eta, \; f^*(F_7(F))=0, \; \omega(F_7(F))=0,\\
[2mm]
f(F_8(F))=0, \; f^*(F_8(F)=f^*(F)(\xi)\eta, \; \omega(F_8(F))=0.
\end{array}$$

\section{The subspaces $h\F$, $v\F$ and $\F_1$ of $\F$}
\emph{The first operator $L_1$.} Let $F\in \F: \;L_1(F)=F-2F_3(F)$.

By straightforward computations, using lemmas 4 and 5, we obtain
\vskip 2mm
{\bf Lemma 6.} $L_1$ is an involutive isometry of $\F$ and commutes with
the action of $U(n)\times 1$.
\vskip 2mm
This lemma and Lemma 1 imply immediately
\vskip 2mm
{\bf Lemma 7.} $\F_1\oplus \F_1^{\perp}$, where
$$\begin{array}{l}
\F_1=L_1^-(\F)=\{F\in \F\,|\,F=F_3(F)\},\\
[2mm]
\F_1^{\perp}=L_1^+(\F)=\{F\in \F\,|\,\omega(F)=0\}.
\end{array}$$
\vskip 2mm
\emph{The second operator $L_2$.} Let $F\in \F_1^{\perp}: \;L_2(F)=F-2\{F_1(F)+F_2(F)\}$.

Analogously to Lemma 6 we obtain
\vskip 2mm
{\bf Lemma 8.} $L_2$ is an involutive isometry of $\F_1^{\perp}$ and commutes
with the action of $U(n)\times 1$. We have
$$\F_1^{\perp}=v\F\oplus h\F \quad ({\rm orthogonally}),$$
where
$$\begin{array}{l}
v\F=L_2^-(\F_1^{\perp})=\{F\in \F\,|\,hF=0, \,\omega(F)=0\},\\
[2mm]
h\F=L_2^+(\F_1^{\perp})=\{F\in \F\,|\,F_1(F)=F_2(F)=0\}.
\end{array}$$
\vskip 2mm
Taking into account lemmas 6, 7 and 8, we obtain a partial decomposition:
\vskip 2mm
{\bf Proposition 1.} $\F=\F_1\oplus v\F\oplus h\F$. The decomposition is orthogonal
and invariant under the action of $U(n)\times 1$. The corresponding components of
$F\in \F$ are
$$p_1(F)=F_3(F), \quad vF=F_1(F)+F_2(F)-2F_3(F), \quad hF.$$

\section{The subspace $v\F$ of $\F$}
\emph{The operator $L_3$.} Let $F\in v\F: \; L_3(F)=F_2(F)-F_1(F)$.
\vskip 2mm
{\bf Lemma 9.} $L_3$ is an involutive isometry of $v\F$ and commutes with
the action of $U(n)\times 1$. We have $v\F=\F_8\oplus (v\F)'$, where
$$\begin{array}{l}
\F_8=L_3^-(v\F)=\{F\in \F\, |\, hF=0, F(x, y, \xi)=0\},\\
[2mm]
(v\F)'=\F_8^{\perp}=L_3^+(v\F)=\{F\in \F\,|\,hF=0, F(\xi,y,z)=0\}.
\end{array}$$
\vskip 2mm
The corresponding components of $F\in v\F$ are $F_1(F)$ and $F_2(F)$.

\emph{The operator $L_4$.} Let $F\in \F_8^{\perp}=(v\F)': \;L_4(F)=-F_4(F)$.
\vskip 2mm
{\bf Lemma 10.} $L_4$ is an involutive isometry of $(v\F)'=\F_8^{\perp}$
and commutes with the action of $U(n)\times 1$. We have
$$(v\F)'=\F_8^{\perp}= \N\F_8^{\perp}\oplus\widetilde{\N}\F_8^{\perp},$$
where
$$\begin{array}{l}
\N\F_8^{\perp}=L_4^-(\F_8^{\perp})=\{F\in \F\,|\,F=F_4(F)\},\\
[2mm]
\widetilde{\N}\F_8^{\perp}=L_4^+(\F_8^{\perp})=\{F\in \F\,|\,F=-F_4(F)\}.
\end{array}$$
\vskip 2mm
The corresponding components  of $F\in (v\F)'=\F_8^{\perp}$ are
$$\frac{1}{2}\{F_2(F)+F_4(F)\}, \quad \frac{1}{2}\{F_2(F)-F_4(F)\}.$$

\emph{The operator $L_5$.} Let $F\in \N\F_8^{\perp}\;(F\in \widetilde \N\F_8^{\perp}):
L_5(F)=-F_5(F)$.
\vskip 2mm
{\bf Lemma 11.} $L_5$ is an involutive isometry of $\N\F_8^{\perp}\;(\widetilde \N\F_8^{\perp})$
and commutes with the action of $U(n)\times 1$. We have
$$\N\F_8^{\perp}=\mathcal Q\mathcal S \F\oplus \mathcal Q \mathcal K \F, \quad
\widetilde \N\F_8^{\perp}=\F_6 \oplus \F_7,$$
where
$$\begin{array}{rl}
\mathcal Q\mathcal S \F &=L_5^-(\N\F_8^{\perp})=\{F\in\F\,|\,F=F_4(F)=F_5(F)\},\\
[2mm]
\mathcal Q \mathcal K \F &=L_5^+(\N\F_8^{\perp})=\{F\in \F\,|\,F=F_4(F)=-F_5(F)\},\\
[2mm]
\F_6&=L_5^-(\widetilde \N\F_8^{\perp})=\{F\in \F\,|\,F=-F_4(F)=F_5(F)\},\\
[2mm]
\F_7&=L_5^+(\widetilde \N\F_8^{\perp})=\{F\in \F\,|\,F=-F_4(F)=-F_5(F)\}.
\end{array}$$

The corresponding components of $F\in \N\F_8^{\perp}\;(F\in \widetilde \N\F_8^{\perp})$
are
$$\begin{array}{l}
\ds{\frac{1}{4}\{F_2(F)+F_4(F)+F_5(F)+F_6(F)\}, \quad
\frac{1}{4}\{F_2(F)+F_4(F)-F_5(F)-F_6(F)\}},\\
[3mm]
\left(\ds{\frac{1}{4}\{F_2(F)-F_4(F)+F_5(F)-F_6(F)\}, \quad
\frac{1}{4}\{F_2(F)-F_4(F)-F_5(F)+F_6(F)\}}\right).
\end{array}$$
\vskip 2mm
\emph{The operator $L_6$.} Let $F\in \mathcal Q\mathcal S \F: \; L_6(F)=F-2F_7(F)$.
\vskip 2mm
{\bf Lemma 12.} $L_6$ is an involutive isometry of $\mathcal Q\mathcal S \F$ and
commutes with the action of $U(n)\times 1$. We have
$$\mathcal Q\mathcal S \F=\F_2\oplus \F_4 \quad ({\rm orthogonally}),$$
where
$$\begin{array}{l}
\F_2=L_6^-(\mathcal Q\mathcal S \F)=\{F\in \F\,|\,F=F_7(F)\},\\
[2mm]
\F_4=L_6^+(\mathcal Q\mathcal S \F)=\{F\in\F\,|\,F=F_4(F)=F_5(F),\;f(F)(\xi)=0\}.\\
\end{array}$$
\vskip 2mm
The corresponding components of $F\in \mathcal Q\mathcal S \F$ in $\F_2$ and
$\F_4$ are
$$F_7(F), \quad \frac{1}{4}\{F_2(F)+F_4(F)+F_5(F)+F_(6)-4F_7(F)\}.$$
\vskip 2mm
\emph{The operator $L_7$.} Let $F\in \mathcal Q\mathcal K \F: \; L_7(F)=F-2F_8(F).$
\vskip 2mm
{\bf Lemma 13.} $L_7$ is an involutive isometry of $\mathcal Q\mathcal K \F$
and commutes with the action of $U(n)\times 1$. We have
$$\mathcal Q\mathcal K \F=\F_3\oplus F_5 \quad (\rm{orthogonally)},$$
where
$$\begin{array}{l}
\F_3=L_7^-(\mathcal Q\mathcal K \F)=\{F\in \F\,|\,F=F_8(F)\},\\
[2mm]
\F_5=L_7^+(\mathcal Q\mathcal K \F)=\{F\in \F\,|\,F=F_4(F)=-F_5(F),\;f^*(F)(\xi)=0\}.
\end{array}$$
\vskip 2mm
The corresponding components of $F\in \mathcal Q\mathcal K \F$ in $\F_3$ and $\F_5$
are
$$F_8(F), \quad \frac{1}{4}\{F_2(F)+F_4(F)-F_5(F)-F_6(F)-4F_8(F)\}.$$
Using lemmas 9\,-\,13, we get
\vskip 2mm
{\bf Proposition 2.} $vF=\F_2\oplus ... \oplus \F_8$. The decomposition is
orthogonal and invariant under the action of $U(n)\times 1$. The corresponding
components of $F\in \F$ in $\F_i\,(i=2,...,8)$ are
$$\begin{array}{l}
p_2(F)=F_7(F),\\
[2mm]
p_3(F)=F_8(F),\\
[2mm]
p_4(F)=\ds{\frac{1}{4}\{F_2(F)+F_4(F)+F_5(F)+F_6(F)-4F_7(F)-F_3(F)\},}\\
[2mm]
p_5(F)=\ds{\frac{1}{4}\{F_2(F)+F_4(F)-F_5(F)-F_6(F)-4F_8(F)-F_3(F)\},}\\
[2mm]
p_6(F)=\ds{\frac{1}{4}\{F_2(F)-F_4(F)+F_5(F)-F_6(F)-F_3(F)\},}\\
[2mm]
p_7(F)=\ds{\frac{1}{4}\{F_2(F)-F_4(F)-F_5(F)+F_6(F)-F_3(F)\}},\\
[2mm]
p_8(F)=F_1(F)-F_3(F).
\end{array}$$

\section{The subspace $h\F$}

Now, let $hV=\{x\in V\,|\,x=hx\}$. Denoting the restrictions of $g$ and $\varphi$
on $hV$ with the same letters, we obtain the Hermitian vector space
$\{hV,g,\varphi\}$ of dimension $2n$. We identify the elements of $h\F$ with
their restrictions on $hV$. Then we can consider the vector space $h\F$ as
the vector space of the tensors $hF$ of type (0,3) over $hV$ having the properties
$$hF(x,y,z)=-hF(x,z,y)=-hF(x, \varphi y, \varphi z)$$
for all $x, y, z \in hV$. The action of $U(n)\times 1$ on $h\F$ coincides
with the action of $U(n)$ on $h\F$. In \cite{GH} the vector space $h\F$ has been
decomposed orthogonally into irreducible components invariant under the action of $U(n)$.

Let $F\in h\F$. It is not difficult to verify that the forms
$$\begin{array}{l}
\begin{array}{rl}
F_9(F)=\ds{\frac{1}{2(n-1)}}&\{g(hx, hy)f(F)(z)-g(hx,hz)f(F)(y)\\
[2mm]
&-g(x, \varphi y)f(F)(\varphi z)+g(x, \varphi z)f(F)(\varphi y)\},\\
\end{array}$$\\
[2mm]
F_{10}(F)=\ds{\frac{1}{2}}\,\{F(x,y,z)+f(\varphi x, \varphi y, z)\},\\
[2mm]
$$\begin{array}{rl}
F_{11}(F)=\ds{\frac{1}{6}}&\{F(x,y,z)+F(y,z,x)+F(z,x,y)\\
[2mm]
&-F(\varphi x, \varphi y,z) -F(\varphi y, \varphi z, x) - F(\varphi z, \varphi x, y)\},
\end{array}$$\\
[2mm]
F_{12}(F)=\ds{\frac{1}{2}}\,\{F(x,y,z)-F(\varphi x, \varphi y,z)\},
\end{array}$$

are also elements of $h\F$.

Using the decomposition in \cite{GH}, we have
\vskip 2mm
{\bf Proposition 3.} $h\F=F_9\oplus F_{10}\oplus F_{11} \oplus F_{11} \oplus F_{12}$,
where
$$\begin{array}{l}
\F_9=\{F\in \F\,|\,F=hF=F_9(F)\},\\
[2mm]
\F_{10}=\{F\in \F\,|\,F=hF=F_{10}(F)-F_9(F)\},\\
[2mm]
\F_{11}=\{F\in \F\,|\,F=hF=F_{11}(F)\},\\
[2mm]
\F_{12}=\{F\in \F\,|\,F=hF=F_{12}(F)-F_{11}(F)\}.
\end{array}$$
The decomposition is orthogonal and invariant under the action of $U(n)\times 1$.
The corresponding components of $F \in \F$ are
$$F_9(F),\quad F_{10}(F)-F_9(F),\quad F_{11}(F), \quad F_{12}(F)-F_{11}(F).$$

\section{Applications to almost contact metric manifolds}

Let $M$ be an almost contact metric manifold with structure $(\varphi, \xi, \eta, g)$,
where $\varphi$ is a tensor field of type (1,1), $\xi$ is a tensor field, $\eta$ is a 1-form,
and $g$ is a Riemannian metric on $M$ such that
$$\begin{array}{l}
\varphi^2x=-x+\eta(x)\,\xi, \quad g(\xi, \xi)=1, \quad \eta \circ \varphi =0,\\
[2mm]
\varphi \,\xi = 0, \quad g(\varphi x, \varphi y)=g(x,y)-\eta(x)\eta(y)
\end{array} $$
for arbitrary vector fields $x, y$ on $M$. For all vector fields $x, y, z $ on $M$ we denote
$$F(x,y,z)=g((\nabla_x\varphi)y,z).\leqno (3)$$

Let $T_pM$ be the tangent space to $M$ at $p\in M$ and $V=T_pM$. The restriction $F_p$
of $F$ given by (3) on $V$ has the properties (1). We shall call $M$ is of class
$\W_i\,(i=1,...,12)$ if $F_p$ is in the subspace $\F_i\,(i=1,...,12)$ for every $p\in M$.
Using the propositions 1, 2 and 3 we obtain 12 basic classes of almost contact metric manifolds .
Further we give the defining conditions for these classes. Let $F$ be given by (3) and
$f$, $f^*$, $\omega$ be $f(F), f^*(F), \omega(F)$ respectively defined by (2).
\vskip 2mm
\hskip 10mm
\emph{The class $\W_1$:}

$$F(x,y,z)=\eta(x)\eta(y)\omega(z)-\eta(x)\eta(z)\omega(y).$$

\vskip 2mm
\hskip 10mm
\emph{The class $\W_2$:}

$$F(x,y,z)=\frac{f(\xi)}{2n}\{\eta(z) g(x,y)-\eta(y)g(x,z)\}.$$
\vskip 1mm
\noindent
This is the class of $\alpha$-Sasakian manifolds.

\vskip 2mm
\hskip 10mm
\emph{The class $\W_3$:}

$$F(x,y,z)=-\frac{f^*(\xi)}{2n}\{\eta(z)g(x,\varphi y)-\eta(y)g(x,\varphi z)\}.$$
\vskip 1mm
\noindent
This is the class of $\alpha$-Kenmotsu manifolds.

\vskip 3mm
\hskip 10mm
\emph{The class $\W_4$:}

$$\begin{array}{rll}
F(x,y,z)&=\eta(y)F(\varphi x, \xi, \varphi z)-\eta(z)F(\varphi x,\xi, \varphi y)&\\
[3mm]
&=\eta(y)F(z,\xi,x)-\eta(z)F(y,\xi,x), & f(\xi)=0.
\end{array}$$

\newpage
\hskip 10mm
\emph{The class $\W_5$:}

$$\begin{array}{rll}
F(x,y,z)&=\eta(y)F(\varphi x, \xi, \varphi z)-\eta(z)F(\varphi x,\xi, \varphi y)&\\
[2mm]
&=-\eta(y)F(z,\xi,x)+\eta(z)F(y,\xi,x), & f^*(\xi)=0.
\end{array}$$

\vskip 1mm
\hskip 10mm
\emph{The class $\W_6$:}

$$\begin{array}{rl}
F(x,y,z)&=-\eta(y)F(\varphi x, \xi, \varphi z)+\eta(z)F(\varphi x,\xi, \varphi y)\\
[2mm]
&=\eta(y)F(z,\xi,x)-\eta(z)F(y,\xi,x).
\end{array}$$

\vskip 1mm
\hskip 10mm
\emph{The class $\W_7$:}

$$\begin{array}{rl}
F(x,y,z)&=-\eta(y)F(\varphi x, \xi, \varphi z)+\eta(z)F(\varphi x,\xi, \varphi y)\\
[2mm]
&=-\eta(y)F(z,\xi,x)+\eta(z)F(y,\xi,x).
\end{array}$$

\vskip 1mm
\hskip 10mm
\emph{The class $\W_8$:}

$$F(hx,hy,hz)=F(x,y,\xi)=0.$$

\vskip 1mm
\hskip 10mm
\emph{The class $\W_9$:}

$$F(\xi,y,z)=F(x,y,\xi)=0,$$
$$\begin{array}{rll}
F(x,y,z)=\ds{\frac{1}{2(n-1)}}&\{g(\varphi x, \varphi y)f(z)
-g(\varphi x, \varphi z)f(y)\\
[2mm]
&-g(x, \varphi y)f(\varphi z)
+g(x, \varphi z)f(\varphi y)\}.
\end{array}$$

\vskip 1mm
\hskip 10mm
\emph{The class $\W_{10}$:}

$$F(\xi,y,z)=F(x,y,\xi)=0,$$
$$F(\varphi x, \varphi y, z)-F(x,y,z)=0, \quad f=0.$$

\vskip 1mm
\hskip 10mm
\emph{The class $\W_{11}$:}

$$F(\xi,y,z)=F(x,y,\xi)=0,$$
$$F(x, x, z)=0.$$

\vskip 1mm
\hskip 10mm
\emph{The class $\W_{12}$:}

$$F(\xi,y,z)=F(x,y,\xi)=0,$$
$$F(x,y,z)+F(y,z,x)+F(z,x,y)=0.$$
\vskip 2mm
The class of cosymplectic manifolds is characterized by $F=0$. This class
is contained in all $\W_i \, (i=1,...,12)$. An almost contact metric manifold
$M$ belongs to two classes $\W_i$, $\W_j$ $(i\neq j)$ iff $M$ is cosymplectic.


\begin{thebibliography}{99}
\bibitem{GH}
Gray A. and L. Hervella. \emph{The Sixteen Classes of Almost Hermitian Manifolds
and Their Linear Invariants.} Ann. Mat. Pura Appl. \textbf{123} (1980), 35-58.
\end{thebibliography}
\end{document}